\documentclass[amstex,12pt, amssymb]{article}

\usepackage{mathtext}
\usepackage[cp1251]{inputenc}
\usepackage[T2A]{fontenc}
\usepackage[dvips]{graphicx}
\usepackage{amsmath}
\usepackage{amssymb}
\usepackage{amsxtra}
\usepackage{latexsym}
\usepackage{ifthen}

\newcounter{lemma}[section]

\newcounter{corollary}[section]

\newcounter{remark}[section]

\newcounter{theorem}[section]

\newcounter{proposition}[section]

\newcounter{example}

\numberwithin{equation}{section}

\begin{document}

\markboth{E~.SEVOST'YANOV, V.~TARGONSKII, D.~ROMASH,
N.~ILKEVYCH}{\centerline{AN ANALOGUE OF KOEBE'S THEOREM IN METRIC
SPACES}}

\def\cc{\setcounter{equation}{0}
\setcounter{figure}{0}\setcounter{table}{0}}

\overfullrule=0pt

%\normalsize\large

\author{EVGENY SEVOST'YANOV\footnote{Corresponding author},
VALERY TARGONSKII, \\ DENYS ROMASH, NATALIYA ILKEVYCH}

\title{
{\bf AN ANALOGUE OF KOEBE'S THEOREM IN METRIC SPACES}}

\date{\today}
\maketitle

%\large
\begin{abstract}
This paper is devoted to the study of mappings in metric spaces. We
investigate mappings satisfying inverse moduli inequalities. We show
that under certain conditions on these mappings, their definition
domains and the spaces in which they act, the image of a ball under
the mappings contains a ball of fixed radius, which corresponds to
the statement of the Koebe theorem on one quarter. As consequences,
we obtain corresponding results in the Sobolev and Orlicz-Sobolev
classes defined in a certain domain of a Riemannian surface or
factor space by the group of fractional-linear mappings of the unit
ball. We also give consequences for manifolds.
\end{abstract}



\section{Introduction}

In our recent papers we have studied problems related to possible
versions of the Koebe-Bloch theorem in the Euclidean $n$-dimensional
space, see, e.g., \cite{ST$_1$}--\cite{ST$_2$}, cf.~\cite{Cr}. We
present, first of all, the classical version of this theorem for
conformal mappings, see, e.g., \cite[Theorem~1.3]{CG}, cf.
\cite{AFW} and \cite{Cr}.

\medskip
{\bf Theorem A.} {\it Let $f:{\Bbb D}\rightarrow {\Bbb C}$ be an
univalent analytic function such that $f(0)=0$ and
$f^{\,\prime}(0)=1.$ Then the image of $f$ covers the open disk
centered at $0$ of radius one-quarter, that is, $f({\Bbb D})\supset
B(0, 1/4).$}

\medskip
In this paper, we aim to present a version of the results
from~\cite{ST$_2$} for metric spaces. In this case, we place the
main emphasis on Loewner spaces, as well as Ahlfors regular spaces
with the Poincare inequality. This ``choice of spaces'' is due to
the fact that it is in them that the modulus of families of paths
and the diameter of the set are metrically related, in other words,
the estimate of the modulus of families of paths joining continua is
expressed through the diameter of these continua (see, for example,
\cite{AS}, \cite{He}). Although in the framework of the manuscript
we prefer to deal with abstract metric spaces, the most important
thing for us is to formulate results in concrete spaces such as
Riemannian surfaces and quotient spaces.

\medskip
Everywhere further $(X, d, \mu)$ and $\left(X^{\,\prime},
d^{\,\prime}, \mu^{\,\prime}\right)$ are metric spaces with metrics
$d$ and $d^{\,\prime}$ and locally finite Borel measures $\mu$ and
$\mu^{\,\prime},$ correspondingly. We will assume that $\mu$ is a
Borel measure such that $0 <\mu(B)<\infty$ for all balls $B$ in $X.$
Let $y_0\in X^{\,\prime},$ $0<r_1<r_2<\infty$ and
\begin{equation}\label{eq1**}
A(y_0, r_1,r_2)=\left\{ y\,\in\,X^{\,\prime}: r_1<d^{\,\prime}(y,
y_0)<r_2\right\}\,,\end{equation}
$$B(y_0, r)=\{y\in X^{\,\prime}: d^{\,\prime}(y,
y_0)<r\}\,, S(y_0, r)=\{y\in X^{\,\prime}: d^{\,\prime}(y,
y_0)=r\}\,.$$
Given sets $E,$ $F\subset X^{\,\prime}$ and a domain $D\subset
X^{\,\prime}$ denote by $\Gamma(E, F, D)$ the family of all paths
$\gamma:[a,b]\rightarrow X^{\,\prime}$ such that $\gamma(a)\in
E,\gamma(b)\in\,F$ and $\gamma(t)\in D$ for $t \in (a, b).$ Given a
domain $D\subset X,$ a {\it mapping} $f:D\rightarrow X^{\,\prime}$
is an arbitrary continuous transformation $x\mapsto f(x).$ Let
$f:D\rightarrow X^{\,\prime},$ let $y_0\in f(D)$ and let
$0<r_1<r_2<d_0=\sup\limits_{y\in f(D)}d^{\,\prime}(y, y_0).$ Now, we
denote by $\Gamma_f(y_0, r_1, r_2)$ the family of all paths $\gamma$
in $D$ such that $f(\gamma)\in \Gamma(S(y_0, r_1), S(y_0, r_2),
A(y_0,r_1,r_2)).$

\medskip
Given a path $\gamma:[a, b]\rightarrow X$ in $(X, d),$ we define its
length as the supremum of the sums
$$
\sum\limits^{k}_{i=1} d(\gamma(t_i),\gamma(t_{i-1}))
$$
over all partitions $a=t_0\leqslant t_1\leqslant\ldots\leqslant
t_k=b$ of the segment $ [a,b].$ A path $\gamma$ is called {\it
rectifiable} if its length is finite. Given a family of paths
$\Gamma$ in $X$, a Borel function $\rho:X\rightarrow[0,\infty]$ is
called {\it admissible} for $\Gamma$, abbr. $\rho\in {\rm
adm}\,\Gamma$, if
\begin{equation*}\label{eq13.2}
\int\limits_{\gamma}\rho\,ds \geqslant 1
\end{equation*}
for all (locally rectifiable) $\gamma\in\Gamma$. Given $p\geqslant
1,$ the $p$-modulus of $\Gamma$ defined as follows:
\begin{equation*}\label{eq13.5}
M_p(\Gamma)=\inf\limits_{\rho\in {\rm
adm}\,\Gamma}\int\limits_X\rho^{\,p}(x)\, d\mu(x)\,.
\end{equation*}
Should ${\rm adm\,}\Gamma$ be empty, we set $M_p(\Gamma)=\infty.$ A
family of paths $\Gamma_1$ in $X$ is said to be {\it minorized} by a
family of paths $\Gamma_2$ in $X,$ abbr. $\Gamma_1>\Gamma_2,$ if,
for every path $\gamma_1\in\Gamma_1$, there is a path
$\gamma_2\in\Gamma_1$ such that $\gamma_2$ is a restriction of
$\gamma_1.$ In this case,
\begin{equation*}\label{eq32*A}
\Gamma_1
> \Gamma_2 \quad \Rightarrow \quad M_p(\Gamma_1)\leqslant M_p(\Gamma_2)
\end{equation*} (see~\cite[Theorem~1]{Fu}).
Let $Q:X^{\,\prime}\rightarrow [0, \infty]$ be a measurable
function. Assume that $D$ and $X^{\,\prime}$ have finite Hausdorff
dimensions $\alpha$ and $\alpha^{\,\prime}\geqslant 1,$
respectively. We will say that {\it $f$ satisfies the inverse
Poletsky inequality} at a point $y_0\in f(D),$ if there is
$r_0=r_0(y_0)>0$ such that the relation
\begin{equation}\label{eq2*A}
M_{\alpha}(\Gamma_f(y_0, r_1, r_2))\leqslant \int\limits_{f(D)\cap
A(y_0,r_1,r_2)} Q(y)\cdot \eta^{{\alpha}^{\prime}}(d^{\,\prime}(y,
y_0))\, d\mu^{\,\prime}(y)
\end{equation}
holds for any $0<r_1<r_2<r_0$ any Lebesgue measurable function
$\eta: (r_1,r_2)\rightarrow [0,\infty ]$ such that
\begin{equation}\label{eqA2}
\int\limits_{r_1}^{r_2}\eta(r)\, dr\geqslant 1\,.
\end{equation}
Note that the inequalities~(\ref{eq2*A}) are well known in the
theory of quasiregular mappings, see e.g.
\cite[Theorem~3.2]{MRV$_1$} or \cite[Theorem~6.7.II]{Ri}.

\medskip
Let $(X, \mu)$ be a metric space with measure $\mu$ and of Hausdorff
dimension $n.$ For each real number $n\geqslant 1,$ we define {\it
the Loewner function} $\Phi_n:(0, \infty)\rightarrow [0, \infty)$ on
$X$ as
\begin{equation}\label{eq2H}
\Phi_n(t)=\inf\{M_n(\Gamma(E, F, X)): \Delta(E, F)\leqslant t\}\,,
\end{equation}
where the infimum is taken over all disjoint nondegenerate continua
$E$ and $F$ in $X$ and
$$\Delta(E, F):=\frac{{\rm dist}\,(E,
F)}{\min\{d(E), d(F)\}}\,.$$
A pathwise connected metric measure space $(X, \mu)$ is said to be a
{\it Loewner space} of exponent $n,$ or an $n$-Loewner space, if the
Loewner function $\Phi_n(t)$ is positive for all $t> 0$ (see
\cite[Section~2.5]{MRSY} or \cite[Ch.~8]{He}). Observe that, ${\Bbb
R}^n$ and ${\Bbb B}^n\subset {\Bbb R}^n$ are Loewner spaces (see
\cite[Theorem~8.2 and Example~8.24(a)]{He}).

\medskip
We agree to say that a metric space $X$ is {\it locally connected}
at a point $x_0$ if for every neighborhood $U$ of $x_0$ there exists
a neighborhood $V$ of the same point such that the set $V\cap X$ is
connected (see \cite[I.6, \S49]{Ku}). We say that a space
$(X,d,\mu)$ is {\it upper $\alpha$-regular at a point} $x_0\in X$ if
there is a constant $C> 0$ such that
$$
\mu(B(x_0,r))\leqslant Cr^{\alpha}$$
for the balls $B(x_0,r)$ centered at $x_0\in X$ with all radii
$r<r_0$ for some $r_0>0.$ We also say that a space  $(X,d,\mu)$ is
{\it upper $\alpha$-regular} if the above condition holds at every
point $x_0\in X.$

\medskip
A space $X^{\,\prime}$ is {\it locally path connected at a point
$x_0\in X^{\,\prime},$} if for an arbitrary neighborhood $U$ of
$x_0$ there exists a neighborhood $V$ of $x_0$ such that $V$ is path
connected. A function $\varphi\colon D\rightarrow{\Bbb R}$
integrable on $B(x_0,r)$ is said to have {\it\,finite mean
oscillation} at a~point $x_0\in D$ (we write $\varphi\in
FMO(x_0)$)~if
$$
\limsup\limits_{\varepsilon\rightarrow
0}\frac{1}{\mu(B(x_0,\varepsilon))} \int\limits_{B(x_0,\varepsilon)}
|\varphi(x)-\overline{\varphi}_{\varepsilon}|\,d\mu(x)<\infty,
$$
where
$$
\overline{\varphi}_{\varepsilon}=\frac{1}{\mu(B(x_0,\varepsilon))}
\int\limits_{B( x_0,\varepsilon)}{\varphi}(x)\,d\mu(x).
$$

\medskip
We say that a domain $ D $ satisfies {\it the condition \textbf{A}}
if any two pairs of (different) points $a\in D^{\,\prime},$
$b\in\overline{D^{\,\prime}}$ and $c\in D^{\,\prime},$
$d\in\overline{D^{\,\prime}}$ may be joined by disjoint paths
$\alpha:[0, 1]\rightarrow \overline{D^{\,\prime}}$ and $\beta:[0,
1]\rightarrow \overline{D^{\,\prime}}$ such that $\alpha(t),
\beta(t)\in D^{\,\prime}$ for any $t\in [0, 1).$ Observe that, the
condition~\textbf{A} always hold for domains in ${\Bbb R}^n,$
locally connected on the boundary, see e.g.
\cite[Proposition~1]{SevSkv}.

\medskip
Given metric spaces $(X, d, \mu)$ and $\left(X^{\,\prime},
d^{\,\prime}, \mu^{\,\prime}\right)$ having Hausdorff dimensions
$\alpha\geqslant 2$ and $\alpha^{\,\prime}\geqslant 2,$
respectively, domains $D\subset X$ and $D^{\,\prime}\subset
X^{\,\prime},$ a continuum $E\subset D,$ $\delta>0$ and a measurable
function $Q:X^{\,\prime}\rightarrow [0, \infty]$ we denote by
$\frak{F}^Q_{E, \delta}(D, D^{\,\prime})$ the family of all open
discrete mappings $f:D\rightarrow D^{\,\prime},$ satisfying
relations~(\ref{eq2*A})--(\ref{eqA2}) at any point $y_0\in
X^{\,\prime}$ such that $d^{\,\prime}(f(E))\geqslant \delta.$ The
following statement holds.

\medskip
\begin{theorem}\label{th1}
{ Let $D$ be a domain in $X,$ while $D$ is locally connected and
locally compact $\alpha$-Loewner space of a Hausdorff dimension
$\alpha\geqslant 2,$ and let $D^{\,\prime}$ be a domain in
$X^{\,\prime}$ satisfying the condition~\textbf{A} which has at
least two different boundary points, while $\overline{D^{\,\prime}}$
is a compactum in $X^{\,\prime}.$ Assume that, $X^{\,\prime}$ is a
locally path connected and locally compact upper
$\alpha^{\,\prime}$-regular metric space of a Hausdorff dimension
$\alpha^{\,\prime}\geqslant 2,$ for which the relation
\begin{equation}\label{eq6} \mu^{\,\prime}(B(z_0, 2r))\leqslant
\gamma\cdot\log^{\alpha^{\,\prime}-2}\frac{1}{r}\cdot
\mu^{\,\prime}(B(z_0, r))
\end{equation}
holds for any $z_0\in X^{\,\prime},$ some $r_0=r_0(z_0)>0$ and every
$r\in (0, r_0).$ Let $\overline{B(x_0, \varepsilon_1)}\subset D$ for
some $\varepsilon_1>0.$ Assume that, $Q\in
FMO(\overline{D^{\,\prime}}).$ Then there is $r_0>0,$ which does not
depend on $f,$ such that
\begin{equation*}\label{eq1D}
f(B(x_0, \varepsilon_1))\supset B(f(x_0), r_0)\qquad \forall\,\,f\in
\frak{F}^Q_{E, \delta}(D, D^{\,\prime})\,,
\end{equation*}
where $B(f(x_0), r_0)=\{w\in X^{\,\prime}: d^{\,\prime}(w,
f(x_0))<r_0\}.$ }
\end{theorem}

\medskip
One can also formulate this theorem for the case when we consider
the mappings in (\ref{eq2*A})--(\ref{eqA2}) with some other
exponents $p$ and $q\geqslant 1,$ generally speaking, not equal to
$\alpha$ and ${\alpha}^{\prime}.$ Given $p, q\geqslant 1,$ we say
that {\it $f$ satisfies the inverse Poletsky inequality with respect
to $(p, q)$-moduli} at a point $y_0\in X^{\,\prime},$ if there is
$r_0=r_0(y_0)>0$ such that the relation
\begin{equation}\label{eq2B}
M_p(\Gamma_f(y_0, r_1, r_2))\leqslant \int\limits_{f(D)\cap
A(y_0,r_1,r_2)} Q(y)\cdot \eta^{q}(d^{\,\prime}(y, y_0))\,
d\mu^{\,\prime}(y)
\end{equation}
holds for any $0<r_1<r_2<r_0$ Lebesgue measurable function $\eta:
(r_1,r_2)\rightarrow [0,\infty ]$ such that~(\ref{eqA2}) holds.

\medskip
A~metric space $(X,d,\mu)$ is called
{\it\,$\widetilde{Q}$-Ahlfors-regular}\index{Ahlfors regular space}
for some $\widetilde{Q}\geqslant 1$ if, for any $x_0\in X$ and some
constant $C\geqslant 1$,
$$
\frac{1}{C}R^{\widetilde{Q}}\leqslant \mu(B(x_0, R))\leqslant
CR^{\widetilde{Q}}\,.
$$
As is well known, Ahlfors $\alpha$-regular spaces have Hausdorff
dimension $\alpha$ (see \cite[p.~61--62]{He}). Let $(X,d,\mu)$ be
a~metric measure space with metric~$d$ and a~locally finite Borel
measure~$\mu$. Following~\cite{He}, \S\,7.22, a~Borel function
$\rho\colon X\rightarrow [0,\infty]$ is said to be an {\it\,upper
gradient} of a~function $u\colon X\rightarrow {\Bbb R}$ if
$$
|u(x)-u(y)|\leqslant \int\limits_{\gamma}\rho\,ds
$$
for any rectifiable path~$\gamma$ joining the points $x$ and $y\in
X$, where, as usual,
$\displaystyle\int\limits_{\gamma}\rho\,ds$~denotes the linear
integral of the~function~$\rho$ over the~path~$\gamma$. We say that
$X$ {\it admits the $(1;p)$-Poincar\'e inequality} if there exist
constants $C\geqslant 1$ and $\tau>0$ such that
$$
\frac{1}{\mu(B)}\int\limits_{B}|u-u_B|\,d\mu(x)\leqslant
C(\operatorname{diam}B) \biggl(\frac{1}{\mu(\tau
B)}\int\limits_{\tau B}\rho^p\,d\mu(x)\biggr)^{1/p}
$$
for any ball $B\subset X$ and arbitrary locally bounded continuous
function $u\colon X\rightarrow {\Bbb R}$ and any upper
gradient~$\rho$ of~$u$, where
$$
u_B:=\frac{1}{\mu(B)}\int\limits_{B}u\,d\mu(x)\,.
$$
Given metric spaces $(X, d, \mu)$ and $\left(X^{\,\prime},
d^{\,\prime}, \mu^{\,\prime}\right)$ having Hausdorff dimensions
$\alpha\geqslant 2$ and $\alpha^{\,\prime}\geqslant 2,$ $p,
q\geqslant 1,$ a domains $D\subset X,$ $D^{\,\prime}\subset
X^{\,\prime},$ a continuum $E\subset D,$ $\delta>0$ and a measurable
function $Q:X^{\,\prime}\rightarrow [0, \infty]$ we denote by
$\frak{F}^{p, q, Q}_{E, \delta}(D, D^{\,\prime})$ the family of open
discrete mappings $f:D\rightarrow D^{\,\prime},$ satisfying
relations~(\ref{eq2*A})--(\ref{eqA2}) at any point $y_0\in
X^{\,\prime}$ such that $d^{\,\prime}(f(E))\geqslant \delta.$ The
following statement holds.

\medskip
\begin{theorem}\label{th1A}
{ Let $D$ be a domain in $X,$ while $D$ is locally connected and
locally compact space of a Hausdorff dimension $\alpha\geqslant 2,$
and let $D^{\,\prime}$ be a domain in $X^{\,\prime},$ while
$\overline{D^{\,\prime}}$ is a compactum in $X^{\,\prime}.$ Assume
that, $X^{\,\prime}$ is a locally path connected and locally compact
upper $\alpha^{\,\prime}$-regular metric space, for which the
relation~(\ref{eq6}) holds for any $z_0\in X^{\,\prime},$ some
$r_0=r_0(z_0)>0$ and every $r\in (0, r_0).$ Let
$\alpha-1<p\leqslant\alpha,$ $0<q\leqslant \alpha^{\,\prime}.$
Assume that, $D$ is Ahlfors regular space supporting $(1;
p)$-Poincar\'{e} inequality. Let $B(x_0, \varepsilon_1)\subset D$
for some $\varepsilon_1>0.$ Assume that, $Q\in
FMO(\overline{D^{\,\prime}})$ and $\frak{F}^{p, q, Q}_{E, \delta}(D,
D^{\,\prime})$ is equicontinuous at any point $x_1\in E.$ Then there
is $r_0>0,$ which does not depend on $f,$ such that
\begin{equation}\label{eq1} f(B(x_0, \varepsilon_1))\supset B(f(x_0), r_0)\qquad
\forall\,\,f\in \frak{F}^{p, q, Q}_{E, \delta}(D, D^{\,\prime})\,,
\end{equation}
where $B(f(x_0), r_0)=\{w\in X^{\,\prime}: d^{\,\prime}(w,
f(x_0))<r_0\}.$ }
\end{theorem}

\section{Preliminaries}

The following result holds (see \cite[Proposition~4.7]{AS}).

\medskip
\begin{proposition}\label{pr_2}
{\, Let $X$ be a $Q$-Ahlfors regular metric measure space that
supports $(1; p)$-Poincar\'{e} inequality for some $p>1$ such that
$Q-1<p\leqslant Q.$ Then there exists a constant $M>0$ having the
property that, for $x\in X,$ $R>0$ and continua $E$ and $F$ in $B(x,
R),$
$$M_p(\Gamma(E, F, X))\geqslant \frac{1}{M}
\cdot\frac{\min\{{\rm diam}\,E, {\rm diam}\,F\}}{R^{1+p-Q}}\,.$$}
\end{proposition}

\medskip
Let $D\subset X,$ $f:D \rightarrow X^{\,\prime}$ be a discrete open
mapping, $\beta: [a,\,b)\rightarrow X^{\,\prime}$ be a path, and
$x\in\,f^{-1}\left(\beta(a)\right).$ A path $\alpha:
[a,\,c)\rightarrow D$ is called a {\it maximal $f$-lifting} of
$\beta$ starting at $x,$ if $(1)\quad \alpha(a)=x\,;$ $(2)\quad
f\circ\alpha=\beta|_{[a,\,c)};$ $(3)$\quad for
$c<c^{\prime}\leqslant b,$ there is no a path $\alpha^{\prime}:
[a,\,c^{\prime})\rightarrow D$ such that
$\alpha=\alpha^{\prime}|_{[a,\,c)}$ and $f\circ
\alpha^{\,\prime}=\beta|_{[a,\,c^{\prime})}.$ The following
statement holds, see e.g.~\cite[Lemma~2.1]{SM}.

\medskip
\begin{lemma}\label{lem9}
{\, Let $X$ and $X^{\,\prime}$ be locally compact metric spaces, let
$X$ be locally connected, let $D$ be a domain in $X,$ and let $f:D
\rightarrow X^{\,\prime}$ be a discrete open mapping. If $\beta:
[a,\,b)\rightarrow X^{\,\prime}$ be a path, and
$x\in\,f^{\,-1}(\beta(a)),$ then there exists a maximal $f$-lifting
of $\beta$ starting at $x.$}
\end{lemma}

\medskip
The following statement was firstly proved in the Euclidean space,
see \cite[Lemma~3.12]{MRV$_2$}. For metric space, it was proved in
some a special form in~\cite[Lemma~2.1]{Skv}, and some latter we
have proved it in the form that we need  (\cite[Lemma~3.1]{FS}).

\medskip
\begin{lemma}\label{lem3}{\,
Let $D$ and $D^{\,\prime}$ be domains with finite Hausdorff
dimensions $\alpha$ and $\alpha^{\,\prime}\geqslant 2$ in spaces
$(X,d,\mu)$ and $(X^{\,\prime},d^{\,\prime}, \mu^{\,\prime}),$
respectively. Assume that $X$ is locally connected and
$\overline{D}, \overline{D^{\,\prime}}$ are compact sets. Let $f$ be
an open discrete mapping of $D$ onto $D^{\,\prime}.$ Let $\beta:[a,
b)\rightarrow D^{\,\prime}$ be a path such that $\beta(t)\rightarrow
\partial D^{\,\prime}$ as $t\rightarrow b-0$ and let
$\alpha:[a,c)\rightarrow D$ be a maximal $f$-lifting of $\beta$
starting at $x\in f^{\,-1}(\beta(a)).$ Then $d(\alpha(t),\partial
D)\rightarrow 0$ as $t\rightarrow c-0.$}
\end{lemma}

\medskip
A space $X$ is called {\it weakly flat} at the point $x_0\in X,$ if,
for any neighborhood $U$ of $x_0$ and for every $P>0,$ there exists
a neighborhood $V\subset U$ of $x_0$ such that
$$M_{\alpha}(\Gamma(E, F, X))\geqslant P$$
for any continua $E, F\subset X$ with $E\cap
\partial U=\varnothing\ne E\cap
\partial V$ and $F\cap \partial U=\varnothing\ne F\cap \partial V.$
A space $X$ is called {\it weakly flat}, if the indicated property
holds for any $x_0\in X.$

\medskip
Given $M>0$ and domains $D\subset X, D^{\,\prime}\subset
X^{\,\prime},$ denote by ${\frak S}_M(D, D^{\,\prime})$ a family of
all open discrete mappings $f$ of $D$ into $D^{\,\prime}$ such that
the conditions~(\ref{eqA2})--(\ref{eq2*A}) hold for any $y_0\in
D^{\,\prime},$ for any $0<r_1<r_2<r_0=r_0(y_0),$ some
$r_0=r_0(y_0)>0,$ some $Q=Q_f$ and $\Vert
Q_f\Vert_{L^1(D^{\,\prime})}\leqslant M.$ The following result is
similar to Theorem~1.1 in~\cite{Skv}; it was proved
in~\cite[Theorem~3.1]{FS} for mappings of $D$ onto $D^{\,\prime},$
which are assumed to be open, discrete, and closed. For mappings of
the class ${\frak S}_M(D, D^{\,\prime})$ mentioned above the proof
is no different and may therefore be omitted.

\medskip
\begin{lemma}\label{lem2}{\,
Let $D$ and $D^{\,\prime}$ be domains with finite Hausdorff
dimensions $\alpha$ and $\alpha^{\,\prime}\geqslant 2$ in spaces
$(X,d,\mu)$ and $(X^{\,\prime},d^{\,\prime}, \mu^{\,\prime}),$
respectively. Assume that $X$ is locally connected, $X^{\,\prime}$
is locally path connected, $\overline{D}$ and
$\overline{D^{\,\prime}}$ are compact sets, and $D$ is weakly flat
as a metric space. Let $D^{\,\prime}\subset X^{\,\prime}$ be a
domain satisfying the condition~\textbf{A} which has at least two
different boundary points. Then the family ${\frak S}_M(D,
D^{\,\prime})$ is equicontinuous in $D.$
 }
\end{lemma}

\medskip
The following statement can be found in~\cite[Lemma~2]{Af},
cf.~\cite[Lemma~13.2]{MRSY}.

\medskip
\begin{proposition}\label{pr3_1}{\,
Let $\left(X^{\,\prime}, d^{\,\prime}, \mu^{\,\prime}\right)$ be a
metric space which is upper $\alpha^{\,\prime}$-regular at $y_0\in
X^{\,\prime}.$ Assume that $Q:X^{\,\prime}\rightarrow [0, \infty]$
belongs to $FMO(y_0).$ If
\begin{equation*}\label{eq7}
\mu^{\,\prime}(B(y_0, 2r))\leqslant
\gamma\cdot\log^{\alpha^{\,\prime}-2}\left(\frac{1}{r}\right)\cdot
\mu^{\,\prime}(B(y_0, r))
\end{equation*}
for some $\gamma, r_0>0$ and every $r\in (0, r_0),$ then
$$\int\limits_{A(y_0, \varepsilon, \varepsilon_0)}\frac{Q(y)\,d\mu^{\,\prime}(y)}
{\left(d(y, y_0)\log\frac{1}{d(y,
y_0)}\right)^{\alpha^{\,\prime}}}=O\left(\log\log\frac{1}{\varepsilon}\right)$$
for some $\varepsilon_0>0.$ }
\end{proposition}

\section{Main Lemmas}

The version of the following lemma was first obtained by the first
two co-authors in~\cite[Lemma~3.2]{ST$_2$}.

\medskip
\begin{lemma}\label{lem1}
{ Assume that, under the conditions of Theorem~\ref{th1}, instead of
condition $Q\in FMO(\overline{D^{\,\prime}}),$ the following
condition is satisfied: for any $y_0\in \overline{D^{\,\prime}}$
there is $\varepsilon_0=\varepsilon_0(y_0)>0$ and a Lebesgue
measurable function $\psi:(0, \varepsilon_0)\rightarrow [0,\infty]$
such that
\begin{equation}\label{eq7***} I(\varepsilon,
\varepsilon_0):=\int\limits_{\varepsilon}^{\varepsilon_0}\psi(t)\,dt
< \infty\quad \forall\,\,\varepsilon\in (0, \varepsilon_0)\,,\quad
I(\varepsilon, \varepsilon_0)\rightarrow
\infty\quad\text{as}\quad\varepsilon\rightarrow 0\,,
\end{equation}
and, in addition,
\begin{equation} \label{eq3.7.2}
\int\limits_{A(y_0, \varepsilon, \varepsilon_0)}
Q(y)\cdot\psi^{\,\alpha^{\,\prime}}(d^{\,\prime}(y,
y_0))\,d\mu^{\,\prime}(y) = o(I^{\alpha^{\,\prime}}(\varepsilon,
\varepsilon_0))\end{equation}
as $\varepsilon\rightarrow 0,$ where $A(y_0, \varepsilon,
\varepsilon_0)$ is defined in (\ref{eq1**}).
Then the conclusion of Theorem~\ref{th1} is true. }
\end{lemma}

\begin{proof}
We primarily use the methodology developed in~\cite{ST$_2$}. Since
by the assumption $D$ is locally compact, we may consider that
$\overline{B(x_0, \varepsilon_1)}$ is compactum in $D.$ Let us prove
the lemma by contradiction. Assume that its conclusion is wrong,
i.e., the relation~(\ref{eq1}) does not hold for any $r_0>0.$ Then
for any $m\in{\Bbb N}$ there is $y_m\in X^{\,\prime}$ and $f_m\in
\frak{F}^Q_{E, \delta}(D, D^{\,\prime})$ such that
$d^{\,\prime}(f_m(x_0), y_m)<1/m$ and $y_m\not\in f_m(B(x_0,
\varepsilon_1)).$ Since $\overline{D^{\,\prime}}$ is a compactum, we
may consider that the sequence $f_m(x_0)$ converges to some $y_0\in
\overline{D^{\,\prime}}$ as $m\rightarrow\infty.$ Now,
$y_m\rightarrow y_0$ as $m\rightarrow\infty,$ as well.

Since by the assumption $d^{\,\prime}(f_m(E))\geqslant \delta$ for
every $m\in{\Bbb N},$ there is $\varepsilon_2>0$ such that
\begin{equation}\label{eq2}
f_m(E)\setminus B(y_0, \varepsilon_2)\ne \varnothing,\qquad
m=1,2,\ldots \,.
\end{equation}
By~(\ref{eq2}), there is $w_m=f_m(z_m)\in X^{\,\prime}\setminus
B(y_0, \varepsilon_2),$ where $z_m\in E.$ Since $E$ is a continuum,
we may consider that $z_m\rightarrow z_0\in E$ as
$m\rightarrow\infty.$ Since $\overline{D^{\,\prime}}$ is a compactum
and the set $\overline{D^{\,\prime}}\setminus B(y_0, \varepsilon_2)$
is closed, we may consider that $w_m\rightarrow w_0\in
\overline{D^{\,\prime}}\setminus B(y_0, \varepsilon_2)$ as
$m\rightarrow\infty.$ Obviously, $w_0\ne y_0.$

\medskip
Since $D$ is Loewner metric space, it is weakly flat (see, e.g.,
\cite[Lemma~2]{SevSkv}). By Lemma~\ref{lem2} the family $f_m$ is
equicontinuous in $D.$ Now, for each $\varepsilon>0$ there is
$\delta=\delta(z_0)>0$ such that $h(f_m(z_0), f_m(z))<\varepsilon$
whenever $d(z, z_0)\leqslant \delta.$ Then, by the triangle
inequality
\begin{equation}\label{eq3}
d^{\,\prime}(f_m(z), w_0)\leqslant d^{\,\prime}(f_m(z), f_m(z_0))+
d^{\,\prime}(f_m(z_0), f_m(z_m))+ d^{\,\prime}(f_m(z_m),
w_0)<3\varepsilon
\end{equation}
for $d(z, z_0)<\delta,$ some $M_1\in {\Bbb N}$ and all $m\geqslant
M_1.$ We may consider that latter holds for any $m=1,2,\ldots .$
Since $w_0\in \overline{D^{\,\prime}}\setminus B(y_0,
\varepsilon_2),$ we may choose $\varepsilon>0$ such that
$\overline{B(w_0, 3\varepsilon)}\cap \overline{B(y_0,
\varepsilon_2)}=\varnothing,$ where $B(w_0, \varepsilon)=\{w\in
X^{\,\prime}: d^{\,\prime}(w, w_0)<\varepsilon\}.$ Then~(\ref{eq3})
implies that
\begin{equation}\label{eq4}
f_m(E_1)\cap \overline{B(y_0, \varepsilon_2)}=\varnothing,\qquad
m=1,2,\ldots \,,
\end{equation}
where $E_1:=\overline{B(z_0, \delta)}.$ Besides that, since
$X^{\,\prime}$ is locally path connected, there is a sequence of
neighborhoods $V_k\subset B(y_0, 2^{\,-k}),$ $k=1,2,\ldots ,$ such
that $V_k$ is path connected. Now, there are increasing subsequences
$m_k,$ $n_k,$ $k=1,2,\ldots,$ of numbers such that $y_{m_k}$ and
$f_{n_k}(x_0)\in V_k.$ Now, $d^{\,\prime}(f_{m_k}(x_0),
y_{m_k})<1/m_k<\frac{1}{k},$ $k=1,2,\ldots .$ Renumbering sequences
$y_{m_k}$ and $f_{n_k}(x_0),$ if required, we may consider that the
latter conditions are satisfied for all $m=1,2,\ldots,$ i.e.,
$y_{m}$ and $f_{m}(x_0)\in V_m,$ $m=1,2,\ldots .$

Join the points $y_m$ and $f_m(x_0)$ by a path $\beta_m:[0,
1]\rightarrow V_m$ such that $\beta_m(0)=f_m(x_0)$ and
$\beta_m(1)=y_m.$ Since $\overline{D}$ is a compact, $D$ is locally
compact, as well. By Lemma~\ref{lem9}, there is a maximal
$f_m$-lifting $\alpha_m:[0, c_m)\rightarrow B(x_0, \varepsilon_1)$
of $\beta_m$ in $B(x_0, \varepsilon_1)$ starting at $x_0.$

\medskip
We prove that $\alpha_m(t)\rightarrow S(x_0, \varepsilon_1)$ as
$t\rightarrow c_m-0,$ cf. \cite[Lemma~3.1]{FS}. First of all, we
show that $S(x_0, \varepsilon_1)\ne\varnothing.$ Indeed, since
$\beta_m(1)=y_m\not \in f_m(B(x_0, \varepsilon_1))$ and
$\beta_m(0)=f_m(x_0)\in f_m(B(x_0, \varepsilon_1)),$ we have that
$|\beta_m|\cap f_m(B(x_0, \varepsilon_1))\ne \varnothing\ne
|\beta_m|\cap(X^{\,\prime}\setminus f_m(B(x_0, \varepsilon_1))).$
Fix $m\in {\Bbb N}.$ Now, by \cite[Theorem~1.I.5.46]{Ku}
$|\beta_m|\cap
\partial f_m(B(x_0, \varepsilon_1))\ne\varnothing.$ Let $\widetilde{y}_0\in
\partial f_m(B(x_0, \varepsilon_1)).$ Then we may find a sequence of
points $y_k\in f_m(B(x_0, \varepsilon_1))$ such that $y_k\rightarrow
\widetilde{y}_0$ as $k\rightarrow\infty.$ Now, there is $x_k\in
B(x_0, \varepsilon_1)$ such that $f_m(x_k)=y_k.$ Since
$\overline{B(x_0, \varepsilon_1)}$ is a compactum, we may also
assume that the sequence $x_k$ converges to some point
$\widetilde{x}_0$ as $k\rightarrow \infty.$ Observe that, the point
$\widetilde{x}_0$ cannot be inner for $B(x_0, \varepsilon_1),$
because otherwise due to the openness of the mapping $f_m$ the point
$\widetilde{y}_0$ is inner for $B(x_0, \varepsilon_1),$ as well.
Then $\widetilde{x}_0\in S(x_0, \varepsilon_1),$ as required.

\medskip
Since we have now established that $S(x_0, \varepsilon_1) \ne
\varnothing,$ only two cases are possible: either $d(\alpha_m(t),
S(x_0, \varepsilon_1))\rightarrow 0$ as $t\rightarrow c_m-0,$ or
$d(\alpha_m(p_k), S(x_0, \varepsilon_1))\geqslant \delta_0> 0$ for
some sequence $p_k\rightarrow c_m-0$ and some $\delta_0>0.$ Let us
show by contradiction that the second situation is impossible. Let
now $d(\alpha_m(p_k), S(x_0, \varepsilon_1))\geqslant \delta_0> 0$
for some sequence $p_k\rightarrow c_m-0$ and some $\delta_0>0.$
Since $\overline{B(x_0, \varepsilon_1)}$ is a compactum, we may
consider that $\alpha(p_k)\rightarrow x_1\in B(x_0, \varepsilon_1)$
as $k\rightarrow\infty.$ Put
$$D_0=\left\{x\in X:  x=\lim_{k\rightarrow \infty} \alpha_m(t_k),\quad
t_k \in [0, c_m), \quad\lim\limits_{k\rightarrow \infty}
t_k=c_m\right\}.$$
Observe that $c_m\ne 1.$  Indeed, in the contrary case
$f_m(x_1)=\beta_m(c)\in f_m(B(x_0, \varepsilon_1))$ that contradicts
the choice of the path $\beta_m.$

\medskip
Let now $c_m\ne 1.$ Passing to the subsequences if necessary, we may
limit ourselves to monotonic sequences $t_k.$ Let $x \in D_0\cap
B(x_0, \varepsilon_1),$ then by the continuity of $f_m$ we obtain
that $f_m(\alpha(t_k))\rightarrow f_m(x) $ as $k\rightarrow \infty,$
where $t_k \in [0, c_m),$ $t_k \rightarrow c_m$ as $k\rightarrow
\infty. $ However, $f_m(\alpha(t_k))=\beta_m(t_k)\rightarrow \beta_m
(c)$ as $k\rightarrow \infty. $ Then the mapping $f_m$ is constant
on $D_0\cap B(x_0, \varepsilon_1).$ Since the sequence of connected
sets $\alpha_m([t_k, c_m))$ is monotone,
$$D_0=\bigcap\limits_{k=1}^{\infty} \overline{\alpha_m([t_k,c))}
\ne\varnothing\,.$$

By~\cite[Theorem~5.II.47.5]{Ku} the set $D_0$ is connected. Let
$L_0$ be a connected component of $D_0\cap B(x_0, \varepsilon_1)$
containing $x_1.$ If $D_0$ contains at least two points, then, by
the definition of connectedness, $L_0$ has the same property. Since
$f_m$ is a discrete mapping and $L_0\subset D,$ the set $L_0$ is
one-point. Hence, $D_0$ is one-point, as well. In this case, the
path $\alpha_m: [0, c_m)\rightarrow B(x_0, \varepsilon_1)$ may be
extended to the closed path $\alpha_m: [0, c_m]\rightarrow B(x_0,
\varepsilon_1),$ and $f_m(\alpha_m(c_m))=\beta_m(c_m).$ Then, by
Lemma~\ref{lem9}, there exists one more maximal $f_m$-lifting
$\alpha_m^{\,\prime}$ of $\beta_m|_{[c_m, 1)}$ starting at the point
$\beta_m(c_m).$ Combining the liftings $\alpha_m$ and
$\alpha^{\,\prime}_m,$ we obtain a new $f_m$-lifting of $\beta_m$
which is defined on $[0, c^{\,\prime}_m),$ $c^{\,\prime}_m\in (c_m,
1).$ This contradicts the maximality of the lifting $\alpha_m.$ The
obtained contradiction indicates that $d(\alpha_m(t), S(x_0,
\varepsilon_1))\rightarrow 0$ as $t\rightarrow c_m-0.$

\medskip
Observe that, $\overline{|\alpha_m|}$ is a continuum in
$\overline{B(x_0, \varepsilon_1)}$ and
$d(\overline{|\alpha_m|})\geqslant d(x_0, S(x_0,
\varepsilon_1))=\varepsilon_1.$ Thus,
\begin{equation}\label{eq1A}
\min\{d(|\alpha_m|), d(E_1)\}\geqslant\min\{\varepsilon_1,
d(E_1)\}\,.
\end{equation}
Besides that, by the triangle inequality, for $y\in |\alpha_m|$ and
$z\in E_1,$ we obtain that
$${\rm dist}\,(|\alpha_m|,
E_1)\leqslant d(y, z)\leqslant d(y, x_0)+d(x_0, z)\leqslant
\varepsilon_1+d(x_0, z)\,.$$
Taking here the $\inf$ over all $z\in E_1,$ we obtain that
\begin{equation}\label{eq2A}
{\rm dist}\,(|\alpha_m|, E_1)\leqslant \varepsilon_1+d(x_0, E_1)\,.
\end{equation}
It follows from~(\ref{eq1A}) and~(\ref{eq2A}) that
\begin{equation}\label{eq3A}
\Delta_m:=\frac{{\rm dist}\,(|\alpha_m|, E_1)}{\min\{d(|\alpha_m|),
d(E_1)\}}\leqslant \frac{\varepsilon_1+d(x_0,
E_1)}{\min\{\varepsilon_1, d(E_1)\}}:=\delta_*>0
\end{equation}
for all $m=1,2,\ldots .$ Let $\Gamma_m:=\Gamma(|\alpha_{m}|, E_1,
D).$ Now, by~(\ref{eq3A}) and due to the definition of the Loewner
function in~(\ref{eq2H})
\begin{equation}\label{eq5}
M_{\alpha}(\Gamma_m)\geqslant \Phi_{\alpha}(\delta_*)>0\,,\qquad
m=1,2,\ldots\,.
\end{equation}
Let us show that the relation~(\ref{eq5}) contradicts the definition
of the mapping $f_m$ in~(\ref{eq2*A})--(\ref{eqA2}). Indeed, recall
that, by the construction $y_{m}$ and $f_{m}(x_0)\in V_m\subset
B(y_0, 2^{\,-m}),$ $m=1,2,\ldots .$ Since $|\beta_m|\in V_m,$
\begin{equation}\label{eq7A}
|\beta_{m}|\subset B(y_0, 2^{\,-m})\,,\qquad m=1,2,\ldots\,.
\end{equation}
Let $m_0\in {\Bbb N}$ be such that $2^{\,-m}<\varepsilon_2,$ where
$\varepsilon_2$ is a number from~(\ref{eq4}). In this case, we
observe that
\begin{equation}\label{eq3G}
f_{m}(\Gamma_{m})>\Gamma(S(y_0, 2^{\,-m}), S(y_0, \varepsilon_2),
A(y_0, 2^{\,-m}, \varepsilon_2))\,.
\end{equation}
Indeed, let $\widetilde{\gamma}\in f_{m}(\Gamma_{m}).$ Then
$\widetilde{\gamma}(t)=f_{m}(\gamma(t)),$ where $\gamma\in
\Gamma_{m},$ $\gamma:[0, 1]\rightarrow D,$ $\gamma(0)\in
|\alpha_{m}|,$ $\gamma(1)\in E_1.$ By the relations~(\ref{eq4}),
and~(\ref{eq7A}), we obtain that $|f_{m}(\gamma(t))|\cap B(y_0,
2^{\,-m})\ne\varnothing \ne |f_{m}(\gamma(t))|\cap
(X^{\,\prime}\setminus B(y_0, 2^{\,-m})).$ Now,
by~\cite[Theorem~1.I.5.46]{Ku} we obtain that, there is $0<t_1<1$
such that $f_{m}(\gamma(t_1))\in S(y_0, 2^{\,-m}).$ Set
$\gamma_1:=\gamma|_{[t_1, 1]}.$ We may consider that
$f_{m}(\gamma(t))\in {\Bbb R}^n\setminus B(y_0, 2^{\,-m})$ for any
$t\geqslant t_1.$ Arguing similarly, we obtain $t_2\in [t_1, 1]$
such that $f_{m}(\gamma(t_2))\in S(y_0, \varepsilon_2).$ Put
$\gamma_2:=\gamma|_{[t_1, t_2]}.$ We may consider that
$f_{m}(\gamma(t))\in B(y_0, \varepsilon_2)$ for any $t\in [t_1,
t_2].$ Now, a path $f_{m}(\gamma_2)$ is a subpath of
$f_{m}(\gamma)=\widetilde{\gamma},$ which belongs to $\Gamma(S(y_0,
2^{\,-m}), S(y_0, \varepsilon_2), A(y_0, 2^{\,-m}, \varepsilon_2)).$
The relation~(\ref{eq3G}) is established.

\medskip
It follows from~(\ref{eq3G}) that
\begin{equation}\label{eq3H}
\Gamma_{m}>\Gamma_{f_{m}}(y_0, 2^{\,-m}, \varepsilon_2)\,.
\end{equation}
Since $I(\varepsilon, \varepsilon_0)\rightarrow\infty$ as
$\varepsilon\rightarrow 0,$ we may consider that $I(2^{\,-m},
\varepsilon_2)>0$ for sufficiently large $m\in {\Bbb N}.$ Set
\begin{equation}\label{eq1_A}\eta_{m}(t)=\left\{
\begin{array}{rr}
\psi(t)/I(2^{\,-m}, \varepsilon_2), & t\in (2^{\,-m}, \varepsilon_2)\,,\\
0,  &  t\not\in (2^{\,-m}, \varepsilon_2)\,,
\end{array}
\right. \end{equation}
where $I(2^{\,-m},
\varepsilon_2)=\int\limits_{2^{\,-m}}^{\varepsilon_2}\,\psi (t)\,
dt.$ Observe that
$\int\limits_{2^{\,-m}}^{\varepsilon_2}\eta_{m}(t)\,dt=1.$ Now, by
the relations~(\ref{eq3.7.2}) and~(\ref{eq3H}), and due to the
definition of $f_{m}$ in~(\ref{eq2*A})--(\ref{eqA2}), we obtain that
$$M_{\alpha}(\Gamma_{m})=M_{\alpha}(\Gamma(|\alpha_{m}|, E_1, D)
)\leqslant M_{\alpha}(\Gamma_{f_{m}}(y_0, 2^{\,-m},
\varepsilon_2))\leqslant$$
\begin{equation}\label{eq3J}
\leqslant \frac{1}{I^{\alpha^{\,\prime}}(2^{\,-m},
\varepsilon_2)}\int\limits_{A(y_0, 2^{\,-m}, \varepsilon_2)}
Q(y)\cdot\psi^{\,\alpha^{\,\prime}}(d^{\,\prime}(y,
y_0))\,d\mu^{\,\prime}(y)\rightarrow 0\quad \text{as}\quad
m\rightarrow\infty\,.
\end{equation}
The relation~(\ref{eq3J}) contradicts with~(\ref{eq5}). The
contradiction obtained above proves the lemma.
\end{proof}

\medskip
{\it Proof of Theorem~\ref{th1}} follows by Lemma~\ref{lem1} and
Proposition~\ref{pr3_1}.~$\Box$

\medskip
The analogue of Lemma~\ref{lem1} in the context of
Theorem~\ref{th1A} is also valid, see below. Since these statements,
as well as their proofs, are very similar, we will present the proof
of the second lemma relatively schematically, emphasizing only the
elements of difference between them.

\medskip
\begin{lemma}\label{lem4}
{ Assume that, under conditions of Theorem~\ref{th1A}, instead of
$Q\in FMO(\overline{D^{\,\prime}}),$ the following condition is
satisfied: for any $y_0\in \overline{D^{\,\prime}}$ there is
$\varepsilon_0=\varepsilon_0(y_0)>0$ and a Lebesgue measurable
function $\psi:(0, \varepsilon_0)\rightarrow [0,\infty]$ such that
the relation~(\ref{eq7***}) holds whenever
\begin{equation} \label{eq3.7.2A}
\int\limits_{A(y_0, \varepsilon, \varepsilon_0)}
Q(y)\cdot\psi^{\,q}(d^{\,\prime}(y, y_0))\,d^{\,\prime}(y) =
o(I^{\,q}(\varepsilon, \varepsilon_0))\end{equation}
as $\varepsilon\rightarrow 0,$ where $A(y_0, \varepsilon,
\varepsilon_0)$ is defined in (\ref{eq1**}).
Then the conclusion of Theorem~\ref{th1} is true. }
\end{lemma}

\begin{proof}
Since by the assumption $D$ is locally compact, we may consider that
$\overline{B(x_0, \varepsilon_1)}$ is compactum in $D.$ Let us prove
the lemma by the contradiction. Assume that its conclusion is wrong,
i.e., the relation~(\ref{eq1}) does not hold for any $r_0>0.$ Then
for any $m\in{\Bbb N}$ there is $y_m\in X^{\,\prime}$ and $f_m\in
\frak{F}^{p, q, Q}_{E, \delta}(D, D^{\,\prime})$ such that
$d^{\,\prime}(f_m(x_0), y_m)<1/m$ and $y_m\not\in f_m(B(x_0,
\varepsilon_1)).$ Since $\overline{D^{\,\prime}}$ is a compactum, we
may consider that the sequence $f_m(x_0)$ converges to some $y_0\in
\overline{D^{\,\prime}}$ as $m\rightarrow\infty.$ Now,
$y_m\rightarrow y_0$ as $m\rightarrow\infty,$ as well.

Since by the assumption $d^{\,\prime}(f_m(E))\geqslant \delta$ holds
for every $m\in{\Bbb N},$ there is $\varepsilon_2>0$ such that
\begin{equation}\label{eq2C}
f_m(E)\setminus B(y_0, \varepsilon_2)\ne \varnothing,\qquad
m=1,2,\ldots \,.
\end{equation}
By~(\ref{eq2C}), there is $w_m=f_m(z_m)\in X^{\,\prime}\setminus
B(y_0, \varepsilon_2),$ where $z_m\in E.$ Since $E$ is a continuum,
we may consider that $z_m\rightarrow z_0\in E$ as
$m\rightarrow\infty.$ Since $\overline{D^{\,\prime}}$ is a compactum
and the set $\overline{D^{\,\prime}}\setminus B(y_0, \varepsilon_2)$
is closed, we may consider that $w_m\rightarrow w_0\in
\overline{D^{\,\prime}}\setminus B(y_0, \varepsilon_2)$ as
$m\rightarrow\infty.$ Obviously, $w_0\ne y_0.$

\medskip
By the assumption, $\frak{F}^{p, q, Q}_{E, \delta}(D, D^{\,\prime})$
is equicontinuous in $E.$ Now, for each $\varepsilon>0$ there is
$\delta=\delta(z_0)>0$ such that $h(f_m(z_0), f_m(z))<\varepsilon$
whenever $d(z, z_0)\leqslant \delta.$ Then, by the triangle
inequality
\begin{equation}\label{eq3B}
d^{\,\prime}(f_m(z), w_0)\leqslant d^{\,\prime}(f_m(z), f_m(z_0))+
d^{\,\prime}(f_m(z_0), f_m(z_m))+ d^{\,\prime}(f_m(z_m),
w_0)<3\varepsilon
\end{equation}
for $d(z, z_0)<\delta,$ some $M_1\in {\Bbb N}$ and all $m\geqslant
M_1.$ We may consider that latter holds for any $m=1,2,\ldots .$
Since $w_0\in \overline{D^{\,\prime}}\setminus B(y_0,
\varepsilon_2),$ we may choose $\varepsilon>0$ such that
$\overline{B(w_0, 3\varepsilon)}\cap \overline{B(y_0,
\varepsilon_2)}=\varnothing,$ where $B(w_0, \varepsilon)=\{w\in
X^{\,\prime}: d^{\,\prime}(w, w_0)<\varepsilon\}.$ Then~(\ref{eq3B})
implies that
\begin{equation*}\label{eq4A}
f_m(E_1)\cap \overline{B(y_0, \varepsilon_2)}=\varnothing,\qquad
m=1,2,\ldots \,,
\end{equation*}
where $E_1:=\overline{B(z_0, \delta)}.$ Besides that, since
$X^{\,\prime}$ is locally path connected, there is a sequence of
neighborhoods $V_k\subset B(y_0, 2^{\,-k}),$ $k=1,2,\ldots ,$ such
that $V_k$ is path connected. Now, there are increasing subsequences
$m_k,$ $n_k,$ $k=1,2,\ldots,$ of numbers such that $y_{m_k}$ and
$f_{n_k}(x_0)\in V_k.$ Now, $d^{\,\prime}(f_{m_k}(x_0),
y_{m_k})<1/m_k<\frac{1}{k},$ $k=1,2,\ldots .$ Renumbering sequences
$y_{m_k}$ and $f_{n_k}(x_0),$ if required, we may consider that the
latter conditions are satisfied for all $m=1,2,\ldots,$ i.e.,
$y_{m}$ and $f_{m}(x_0)\in V_m,$ $m=1,2,\ldots .$

Join the points $y_m$ and $f_m(x_0)$ by a path $\beta_m:[0,
1]\rightarrow V_m$ such that $\beta_m(0)=f_m(x_0)$ and
$\beta_m(1)=y_m.$ Since $\overline{D}$ is a compact, $D$ is locally
compact, as well. By Lemma~\ref{lem9}, there is a maximal
$f_m$-lifting $\alpha_m:[0, c_m)\rightarrow B(x_0, \varepsilon_1)$
of $\beta_m$ in $B(x_0, \varepsilon_1)$ starting at $x_0.$ Reasoning
similarly to the proof of Lemma~\ref{lem1}, we may prove that
$\alpha_m(t)\rightarrow S(x_0, \varepsilon_1)$ as $t\rightarrow
c_m-0.$ Observe that, $\overline{|\alpha_m|}$ is a continuum in
$\overline{B(x_0, \varepsilon_1)}$ and
$d(\overline{|\alpha_m|})\geqslant d(x_0, S(x_0,
\varepsilon_1))=\varepsilon_1.$ Thus,
\begin{equation}\label{eq1B}
\min\{d(|\alpha_m|), d(E_1)\}\geqslant\min\{\varepsilon_1,
d(E_1)\}\,.
\end{equation}
Since $|\alpha_m|\subset \overline{B(x_0, \varepsilon_1)}$ and $E$
is a continuum in $D,$ there is $R>0$ such that $|\alpha_m|\cup
E\subset B(x_0, R).$ Let $\Gamma_m:=\Gamma(|\alpha_{m}|, E_1, D).$
Now, by Proposition~\ref{pr_2} and~(\ref{eq1B})
\begin{equation}\label{eq5A}
M_{p}(\Gamma_m)\geqslant \frac{1}{M} \cdot\frac{\min\{\varepsilon_1,
d(E_1)\}}{R^{1+p-\alpha}}>0\,,\qquad m=1,2,\ldots\,,
\end{equation}
where $M>0$ is some constant. Let us show that the
relation~(\ref{eq5A}) contradicts the definition of the mapping
$f_m$ in~(\ref{eq2*A})--(\ref{eqA2}). Arguing similarly
to~(\ref{eq7A}) and~(\ref{eq3G}), we obtain that
\begin{equation}\label{eq3G_1}
f_m(\Gamma_{m})>\Gamma(S(y_0, 2^{\,-m}), S(y_0, \varepsilon_2),
A(y_0, 2^{\,-m}, \varepsilon_2))\,.
\end{equation}
It follows from~(\ref{eq3G_1}) that
\begin{equation}\label{eq3H_1}
\Gamma_{m}>\Gamma_{f_{m}}(y_0, 2^{\,-m}, \varepsilon_2)\,.
\end{equation}
Choosing the function $\eta_m$ similarly to~(\ref{eq1_A}) and using
the relations~(\ref{eq3.7.2A}) and~(\ref{eq3H_1}) and the definition
of $f_{m}$ in~(\ref{eq2*A})--(\ref{eqA2}), we obtain that
\begin{gather}\nonumber M_{p}(\Gamma_m)=M_{p}(\Gamma(|\alpha_{m}|, E_1, D) )\leqslant
M_{p}(\Gamma_{f_{m}}(y_0, 2^{\,-m}, \varepsilon_2))\\
\label{eq3J_1} \leqslant \frac{1}{I^{q}(2^{\,-m},
\varepsilon_2)}\int\limits_{A(y_0, 2^{\,-m}, \varepsilon_2)}
Q(y)\cdot\psi^{\,q}(d^{\,\prime}(y,
y_0))\,d\mu^{\,\prime}(y)\rightarrow 0\quad \text{as}\quad
m\rightarrow\infty\,.
\end{gather}
The relation~(\ref{eq3J_1}) contradicts with~(\ref{eq5A}). The
contradiction obtained above proves the lemma.
\end{proof}

\medskip
{\it Proof of Theorem~\ref{th1A}} follows by Lemma~\ref{lem2} and
Proposition~\ref{pr3_1}.

Indeed, we set
$\psi(t)=\left(\frac{1}{t\log\frac{1}{t}}\right)^{\frac{\alpha^{\,\prime}}{q}}.$
Since by the assumption $q\leqslant \alpha^{\,\prime},$ we have that
$\psi(t)\geqslant \frac {1}{t\,\log{\frac1t}}$ for sufficiently
small $t>0$ and, consequently,
$$I(\varepsilon,
\varepsilon_0)\,:=\,\int\limits_{\varepsilon}^{\varepsilon_0}\psi(t)\,dt\,\geqslant
\log{\frac{\log{\frac{1}
{\varepsilon}}}{\log{\frac{1}{\varepsilon_0}}}}\,.$$
Now, by the latter relation and by Proposition~\ref{pr3_1} we obtain
that
$$\int\limits_{A(y_0, \varepsilon, \varepsilon_0)}
Q(y)\cdot\psi^{\,q}(d^{\,\prime}(y,
y_0))\,d\mu^{\,\prime}(y)\leqslant C_1\cdot
\log\log\frac{1}{\varepsilon}$$
for some $C_1>0$ and sufficiently small $\varepsilon>0.$ The latter
implies the validity of the relation~(\ref{eq3.7.2A}). Thus, the
desired conclusion follows by Lemma~\ref{lem4}.~$\Box$

\section{Convergence theorems}

The following statement holds, cf. \cite{ST$_1$}--\cite{ST$_2$} and
\cite{Cr}.

\medskip
\begin{theorem}\label{th2}
{\,Assume that, all conditions of Theorem~\ref{th1}
(Theorem~\ref{th1A}) are met. Let $f_j:D\rightarrow X^{\,\prime},$
$j=1,2,\ldots,$ be a sequence of open discrete mappings satisfying
the conditions~(\ref{eq2*A})--(\ref{eqA2}) (respectively,
relations~(\ref{eq2B}), (\ref{eqA2})) at any point $y_0\in
X^{\,\prime}$ and converging to some mapping $f:D\rightarrow
X^{\,\prime}$ as $j\rightarrow\infty$ locally uniformly in $D.$ Then
either $f$ is a constant, or $f$ is light and open mapping
$f:D\rightarrow X^{\,\prime}.$}
\end{theorem}

\begin{proof}
Assume that $f$ is not a constant. Now, the lightness of $f$ under
conditions of Theorem~\ref{th3} follows by~\cite[Theorem~1]{SevSkv}.
The lightness of $f$ under the conditions of Theorem~\ref{th1} may
be given similarly to~\cite[Theorem~1]{SevSkv}.

\medskip
It remains to show that $f$ is open. We use now the scheme of the
proof from~\cite{ST$_2$}. Let $A$ be an open set and let $x_0\in A.$
We need to show that, there is $\varepsilon^*>0$ such that
$B(f(x_0), \varepsilon^*)\subset f(A).$ Since $A$ is open, there is
$\varepsilon_1>0$ such that $\overline{B(x_0, \varepsilon_1)}\subset
A$ for some $\varepsilon_1>0.$ Since $D$ is locally compact, we may
consider that $\overline{B(x_0, \varepsilon_1)}$ is a compactum in
$D.$ By Theorem~\ref{th1} (respectively, Theorem~\ref{th1A})
\begin{equation}\label{eq1AA}
f_j(B(x_0, \varepsilon_1))\supset B(f_j(x_0), r_0)
\end{equation}
for any $j\in {\Bbb N}$ and some $r_0>0.$ Let $y\in B(f(x_0),
r_0/2).$ By~(\ref{eq1AA}), $y=f_m(x_m)$ for some $x_m\in B(x_0,
\varepsilon_1).$ Due to the compactness of $\overline{B(x_0,
\varepsilon_1)},$ we may consider that $x_m\rightarrow z_0\in
\overline{B(x_0, \varepsilon_1)}$ as $m\rightarrow\infty.$ By the
continuity of $f$ in $A,$ since $\overline{B(x_0,
\varepsilon_1)}\subset A,$ we obtain that $f(x_m)\rightarrow f(z_0)$
as $m\rightarrow\infty.$ So, we have that $f(x_m)\rightarrow f(z_0)$
as $m\rightarrow\infty$ and simultaneously $y=f_m(x_m)$ for
sufficiently large $m\in {\Bbb N}.$ Thus, by the triangle
inequality,
$$d^{\,\prime}(y, f(z_0))=d^{\,\prime}(f_m(x_m), f(z_0))\leqslant$$
$$\leqslant d^{\,\prime} (f_m(x_m), f(x_m))+d^{\,\prime} (f(x_m), f(z_0))\rightarrow 0\,,$$
$m\rightarrow\infty.$ Thus, $y=f(z_0)\in f(\overline{B(x_0,
\varepsilon_1)})\subset f(A).$ So, $y\in f(A),$ i.e., $B(f(x_0),
r_0/2)\subset f(A),$ as required. \end{proof}

\section{Some applications for Sobolev and Orlicz-Sobolev classes on factor spaces}

This section is devoted to the study of mappings on quotient spaces
corresponding to a certain discontinuous group of M\"{o}bius
transformations of the unit ball. All necessary definitions of this
section, including the definitions of factor-spaces, metrics and
measures on it, can be found, for instance, in~\cite{RV} and
\cite{Sev$_2$}.

\medskip
The following definitions and notions are from~\cite{Sev$_2$}. Let
$G$ be some group of M\"{o}bius transformations of the unit ball
${\Bbb B}^n$ onto itself. In what follows, the points $x$ and
$y\in{\Bbb B}^n$ will be called {\it $G$-equivalent} (or shorter,
{\it equivalent}), if there is $A\in G$ such that $x=A(y).$ A set
consisting of equivalence classes of elements according to the
indicated principle is denoted by ${\Bbb B}^n/G.$ Denote by
${\mathcal{G M}}({\Bbb B}^n)$ the group of all M\"{o}bius
transformations of ${\Bbb B}^n$ onto itself. According
to~\cite[Section~3.4]{MS}, the {\it hyperbolic measure} of the
Lebesgue measurable set $A\subset{\Bbb B}^n$ is determined by the
relation
\begin{equation}\label{eq2B_1}
v(A)=\int\limits_A\frac{2^n\, dm(x)}{{(1-|x|^2)}^n}\,.
\end{equation}
We define the {\it hyperbolic distance} $h(x, y)$ between the points
$x, y\in{\Bbb B}^n$ by the relation
\begin{equation}\label{eq3_1}
h(x, y)=\log\,\frac{1+t}{1-t}\,,\quad
t=\frac{|x-y|}{\sqrt{|x-y|^2+(1-|x|^2)(1-|y|^2)}}\,,
\end{equation}
see, for example, \cite[relation~(2.18), Remark~2.12 and
Exercise~2.52]{Vu}. Note that $h(x, y)=h(g(x), g(y))$ for any $g\in
{\mathcal{G M}}({\Bbb B}^n),$ see~\cite[relation~(2.20)]{Vu}.
In what follows, we denote by $I$ the identity mapping in~${\Bbb
R}^n.$ According to~\cite[Section~3.4]{MS}, the set of the form
\begin{equation}\label{eq1_1}
P=\{x\in {\Bbb B}^n: h(x, x_0)<h(x, T(x_0))\quad {\rm
for\,\,all}\quad T\in G\setminus\{I\}\}
\end{equation}
is called a {\it normal fundamental polyhedron} with center at the
point $x_0.$ Let $\pi:{\Bbb B}^n\rightarrow {\Bbb B}^n/G$ be the
natural projection of ${\Bbb B}^n$ onto the factor space ${\Bbb
B}^n/G,$ then the {\it hyperbolic measure} of the set $A\subset
{\Bbb B}^n/G$ is defined by the relation $v(P\cap\pi^{\,-1}(A)),$
where $P$ is a normal fundamental polyhedron~(\ref{eq1_1}). It is
easy to verify by direct calculations that the hyperbolic measure
does not change under any map~$g\in{\mathcal{G M}}({\Bbb B}^n).$

\medskip
Every point $p\in {\Bbb B}^n/G$ can be identified with a set
$$p=G_z=\{\xi\in {\Bbb B}^n:\,\xi=g(z)\quad{\rm for\,\,some}\quad g\in G
\}\,,$$
called the {\it orbit} of the point $z\in\pi^{\,-1}(p).$ This
definition is independent of the choice of $z\in\pi^{\,-1}(p)$ and
therefore for given elements
$$p_1=G_{z_1}\in {\Bbb B}^n/G\quad {\rm and}\quad p_2=G_{z_2}\in {\Bbb
B}^n/G$$
we may set
\begin{equation}\label{eq2_1}
\widetilde{h}(p_1, p_2):=\inf\limits_{g_1, g_2\in G}h(g_1(z_1),
g_2(z_2))\,.
\end{equation}
The length of the path $\gamma:[a, b]\rightarrow {\Bbb B}^n/G$ on
the segment $[a, t],$ $a\leqslant t\leqslant b,$ is defined as
follows:
\begin{equation*}\label{eq5D}
l_{\gamma}(t):=\sup\limits_{\pi}\sum\limits_{k=0}^{m-1}\widetilde{h}(\gamma(t_k),
\gamma(t_{k+1}))\,,
\end{equation*}
where~$\sup$ is taken over all subdivisions $\pi=\{a=t_0\leqslant
t_1\leqslant t_2\leqslant\ldots\leqslant t_m=t\}.$
Let $G$ be a group of M\"{o}bius automorphisms of the unit ball. We
say that $G$ acts {\it discontinuously} on ${\Bbb B}^n,$ if each
point $x\in {\Bbb B}^n$ has a neighborhood $U$ such that $g(U)\cap
U=\varnothing$ for all $g\in G,$ $g\ne I,$ except maybe a finite
number of elements $g.$ We say that $G$ does not have {\it fixed
points} in $ {\Bbb B}^n,$ if for any $a\in{\Bbb B}^n$ the equality
$g(a)=a$ is possible only when $g=I.$

\medskip
Let $D$ be a domain in ${\Bbb B}^n/G$ and let $Q:{\Bbb B}^n/
G\rightarrow [0, \infty]$ be a measurable function with a respect to
measure $\widetilde{v}.$ Due to~\cite[Section~7]{MRSY}, a mapping
$f:D\rightarrow D_*$ is called a {\it ring $Q$-mapping at the point
$p_0\in \overline{D}$}, if the relation
\begin{equation}\label{eq1F}
M(f(\Gamma(S_1, S_2, D)))\leqslant \int\limits_{\widetilde{A}}
Q(p)\cdot \eta^n(\widetilde{h}(p, p_0))\,d\widetilde{v}(p)
\end{equation}
holds for some $0<r_0<\infty$ and any $0<r_1<r_2<r_0,$ where
$S_1=\widetilde{S}(p_0, r_1)$, $S_2=\widetilde{S}(p_0, r_2)$,
$\widetilde{A}=\widetilde{A}(p_0, r_1, r_2)=\{p\in {\Bbb B}^n/G:
r_1<\widetilde{h}(p, p_0)<r_2\}$ and $\eta:(r_1, r_2)\rightarrow [0,
\infty]$ is arbitrary Lebesgue measurable function with
\begin{equation}\label{eq*3!!}
\int\limits_{r_1}^{r_2}\eta(r)\,dr\geqslant 1\,.
\end{equation}
Let $D$ and $D_{\,*}$ be domains of ${\Bbb B}^n/G$ and ${\Bbb
B}^n/G_{\,*},$ respectively. Suppose that ${\Bbb B}^n/G$ and ${\Bbb
B}^n/G_{\,*}$ are metric spaces with metrics $\widetilde{h}$ and
$\widetilde{h}_*.$ Hereinafter, $\widetilde{h}$ and
$\widetilde{h}_*$ are determined solely by relation~(\ref{eq2_1}).
(It will be established below, under what conditions on $G$ and
$G_*$ the indicated functions $\widetilde{h}$ and $\widetilde{h}_*$
are really metrics). Let $ds_{\widetilde{h}}$ and $d\widetilde{v}$
be elements of length and volume on ${\Bbb B}^n/G,$ and let
$ds_{\widetilde{h}_*}$ and $d\widetilde{v}_*$ be elements of length
and volume on ${\Bbb B}^n/G_*,$ correspondingly. If $n=2,$ we call
${\Bbb B}^n/G$ and ${\Bbb B}^n/G_{\,*}$ {\it Riemannian surfaces}
and denote ${\Bbb B}^2/G={\Bbb S}$ and ${\Bbb B}^2/G_{\,*}={\Bbb
S}_{\,*}.$ The following statements hold.

\medskip
\begin{proposition}\label{pr1A}{\it\,{\rm (\cite[Proposition~1.1]{Sev$_2$})}.\,
Suppose that $G$ is a group of M\"{o}bius transformations of the
unit ball ${\Bbb B}^n,$ $n\geqslant 2,$ onto itself, acting
discontinuously and not having fixed points in ${\Bbb B}^n.$ Then
the factor space ${\Bbb B}^n/G$ is a conformal manifold, that is, a
topological manifold in which any two charts are interconnected by
means of conformal mappings. At the same time, the natural
projection $\pi,$ which maps ${\Bbb B}^n$ onto ${\Bbb B}^n/G,$ is a
local homeomorphism. Moreover, the corresponding pairs of the form
$(U, \pi^{\,-1}),$ where $U$ is some neighborhood of an arbitrary
point $p\in {\Bbb B}^n/G,$ in which the mapping $\pi^{\,-1}$ is
well-defined and continuous, can be considered as charts
corresponding to the specified manifold.}
\end{proposition}

\medskip
\begin{proposition}\label{pr5}
{\,{\rm (\cite[Lemma~2.2]{Sev$_2$})}. Suppose that $G$ is a group of
M\"{o}bius transformations of the unit ball ${\Bbb B}^n,$
$n\geqslant 2,$ onto itself, acting discontinuously and not having
fixed points in ${\Bbb B}^n.$ Then the space~${\Bbb B}^n/G$ is
metrizable, and the corresponding metric can be determined by
relation~(\ref{eq2_1}).}
\end{proposition}

\medskip
Given $z_1, z_2\in {\Bbb B}^n,$ we set
\begin{equation*}
d(z_1, z_2):=\widetilde{h}(\pi(z_1), \pi(z_2))\,,
\end{equation*}
where $\widetilde{h}$ is defined in~(\ref{eq2_1}). Note that, by the
definition, $d(z_1, z_2)\leqslant h(z_1, z_2).$ The following
statement is true.

\medskip
\begin{proposition}\label{pr4}
{\, {\rm (\cite[Lemma~2.3]{Sev$_2$})}. Suppose that $G$ is a group
of M\"{o}bius transformations of the unit ball ${\Bbb B}^n,$
$n\geqslant 2,$ onto itself, acting discontinuously and not having
fixed points in ${\Bbb B}^n.$ Then for any compact set $A\subset
{\Bbb B}^n$ there is $\delta=\delta(A)>0$ such that
\begin{equation*}\label{eq34}
d(z_1, z_2)=h(z_1, z_2) \quad{\rm for\,\, all}\quad z_1, z_2\in
A\quad {\rm such\,\, that}\quad h(z_1, z_2)<\delta\,.
\end{equation*}
}
\end{proposition}

\medskip
\begin{proposition}\label{pr5_1}
{\, {\rm (\cite[Lemma~2.4]{Sev$_2$})}. Let $0<2r_0<1,$ then there is
$C_1=C_1(r_0)>0$ such that
\begin{equation}\label{eq9B}
C_1\cdot h(z_1, z_2)\leqslant |z_1-z_2|\leqslant  h(z_1, z_2)\quad
{\rm for\,\,all}\quad z_1, z_2\in B(0, r_0)\,. \end{equation}
Moreover, the right-hand inequality in~(\ref{eq9B}) holds for all
$z_1, z_2\in {\Bbb B}^n.$
}
\end{proposition}

\medskip
Let us to prove the following statement.

\medskip
\begin{proposition}\label{pr3}{\,
Suppose that $G$ is a group of M\"{o}bius transformations of the
unit ball ${\Bbb B}^n,$ $n\geqslant 2,$ onto itself, acting
discontinuously and not having fixed points in ${\Bbb B}^n.$ Then
the factor space ${\Bbb B}^n/G$ is locally Ahlfors regular, i.e.,
for any $p_0\in {\Bbb B}^n/G$ there is $C>0$ and $r_0>0$ such that
\begin{equation}\label{eq1C}
\frac{1}{C}r^n\leqslant \widetilde{v}(B_{\widetilde{h}}(p_0,
r))\leqslant Cr^n
\end{equation}
for any $r\in (0, r_0).$}
\end{proposition}

\medskip
\begin{proof}
Let $p_0\in {\Bbb B}^n/G$ and $z_0\in {\Bbb B}^n$ be such that
$\pi(z_0)=p_0,$ where $\pi$ is the natural projection of ${\Bbb
B}^n$ onto ${\Bbb B}^n/G.$ By Propositions~\ref{pr1A}, \ref{pr4}
and~\ref{pr5_1} we may choose $\varepsilon_0>0$ such that $\pi$ maps
$B_h(z_0,\varepsilon_0)$ onto~$\widetilde{B}(p_0, \varepsilon_0)$
homeomorphically and isometrically. Without loss of generality, we
may assume that $z_0=0,$ see comments before Remark~2.1
in~\cite{Sev$_2$}. Now,
$\widetilde{v}(\pi^{\,-1}(B_{\widetilde{h}}(p_0, r)))=v(B_h(0, r)),$
where $v$ is defined in~(\ref{eq2B_1}). Observe that, by
Proposition~\ref{pr5_1} $(B(0, C_1r))\subset B_h(0, r)\subset B(0,
r),$ where $B(0, r)$ denotes the Euclidean ball of the radius $r$
centered at the origin, and $B_h(0, r)$ is a hyperbolic ball,
$B_h(0, r)=\{x\in {\Bbb D}: h(x, 0)<r\},$ corresponding to the
metric $h$ defined in~(\ref{eq3_1}). Let $0<r<r_0<1.$ Now,
by~(\ref{eq2B_1}) we obtain that
$$\widetilde{v}(B_{\widetilde{h}}(p_0, r))=v(B_h(0, r))=
\int\limits_{B_h(0, r)}\frac{2^n\, dm(x)}{{(1-|x|^2)}^n}\leqslant
\frac{2^n\Omega_nr^n}{(1-r_0^2)^n}$$
and
$$\widetilde{v}(B_{\widetilde{h}}(p_0, r))=v(B_h(0, r))=
\int\limits_{B_h(0, r)}\frac{2^n\, dm(x)}{{(1-|x|^2)}^n}\geqslant
\int\limits_{B_h(0, r)} 2^n\, dm(x)\geqslant
2^n\Omega_nC^n_1r^n\,,$$
where $C_1$ is a constant from Proposition~\ref{pr5_1}. Now, we have
proved~(\ref{eq1C}) with $$C=\max\biggl\{2^n\Omega_nC^n_1,
\frac{2^n\Omega_n}{(1-r_0^2)^n}\biggr\}\,.$$
\end{proof}

\medskip
Let $D$ and $D_{\,*}$ be domains on the factor spaces ${\Bbb B}^n/G$
and ${\Bbb B}^n/G_{\,*},$ respectively, and let $f:D\rightarrow
D_{\,*}.$ We say that $f\in W_{\rm loc}^{1,1}(D),$ if for each point
$x_0\in D$ there are open neighborhoods $U$ and $V,$ containing the
points $x_0$ and $f(x_0),$ respectively, in which the natural
projections $\pi:\pi^{-1}(U)\rightarrow U$ and
$\pi_*:\pi_*^{-1}(V)\rightarrow V$ are one-to-one mappings, while
$F=\pi_*^{\,-1}\circ f\circ\pi\in W_{\rm loc}^{1,1}(\pi^{\,-1}(U)).$

\medskip
We write $f\in W_{\rm loc}^{1,p}(D),$ $p\geqslant 1,$ if $f\in
W_{\rm loc}^{1,1}(D)$ and, in addition, $\frac{\partial
f_i}{\partial x_j}\in L^p_{\rm loc}(D)$ in local coordinates. For a
given mapping $f:D\rightarrow{\Bbb R}^n,$ which is differentiable
almost everywhere in $D,$ we set
$$\Vert
f^{\,\prime}(x)\Vert=\max\limits_{|h|=1}{|f^{\,\prime}(x)h|}\,,\quad
J(x, f)={\rm det}\, f^{\,\prime}(x)\,.$$
The {\it outher dilatation} $K_O(x, f)$ of the mapping $f$ at the
point $x$ is defined by the relation
\begin{equation}\label{eq0.1.1}
K_O(x,f)\quad =\quad\left\{
\begin{array}{rr}
\frac{{\Vert
f^{\,\prime}(x)\Vert}^n}{|J(x,f)|}, & J(x,f)\ne 0,\\
1,  &  f^{\,\prime}(x)=0, \\
\infty, & \text{otherwise}
\end{array}
\right.\,.
\end{equation}
If we are talking about the mapping $f$ of domains $D$ and $D_*,$
belonging to the factor-spaces  ${\Bbb B}^n/G$ and ${\Bbb B}^n/G_*,$
respectively, then we set $K_O(p, f)=K_O(\varphi(p), F),$ where
$F=\psi\circ f\circ\varphi^{-\,1},$ $(U, \varphi)$ are local
coordinates of $p$ and $(V, \psi)$ are local coordinates of~$f(p).$
By virtue of Proposition~\ref{pr1A}, the mappings $\pi$ and $\pi_*$
can be considered as such local coordinates, besides that, $K_O(p,
f)$ does not depend on the choice of local coordinates, because the
outher dilatation of the conformal mapping is equal to one. The
following results holds, see \cite[Lemma~3.1]{RV}.

\medskip
\begin{proposition}\label{pr1}
{\, Let $D$ and $D_*$ be domains on Riemannian surfaces ${\Bbb S}$
and ${\Bbb S}_*.$ If $f:D\rightarrow D_*$ is a homeomorphism of
finite distortion with $K_O(p, f)\in L^1_{\rm loc},$ then $f$
satisfies~(\ref{eq1F})--(\ref{eq*3!!}) for $n=2$ and $Q(p)=K_O(p,
f).$}
\end{proposition}

\medskip
Let $D$ and $D_*$ be domains in ${\Bbb B}^n/G$ and ${\Bbb B}^n/G_*$,
correspondingly, and let $\varphi\colon
[0,\infty)\rightarrow[0,\infty)$ be a non-decreasing function. We
say that $f\in W^{1,\varphi}_{\rm loc}(D)$, if for any $p\in D$ and
$f(p)$ there are neighborhoods $U$ and $V$ of these points and the
coordinate mappings $\pi:(\pi^{\,-1}(U))\rightarrow U$ and
$\pi_*:(\pi_*^{\,-1}(V))\rightarrow V$, which are one-to-one in $U$
and $V$, respectively, such that $F=\pi_*^{\,-1}\circ f\circ\pi\in
W_{\rm loc}^{1,1}(\pi^{\,-1}(U))$, while
\begin{equation*}\label{eq1AAA}
\int_{\pi^{\,-1}(U)}\varphi\left(|\nabla
F(x)|\right)\,dm(x)<\infty\,,
\end{equation*}
where, as usual,
$$|\nabla
F(x)|=\sqrt{\sum\limits_{i=1}^n\sum\limits_{j=1}^n\left(\frac{\partial
F_i}{\partial x_j}\right)^2}\,.$$
The following results holds, see \cite[Theorem~2]{Sev$_1$}.

\medskip
\begin{proposition}\label{pr2}
{\, Let $n\geqslant 3$, and let $G$ be a some group of M\"{o}bius
automorphisms of the unit ball which acts discontinuously and have
no fixed points in ${\Bbb B}^n.$ Assume that, the group $G_*$, which
corresponds to the space ${\Bbb B}^n/G_*$, consists only from one
mapping $g(x)=x$. Let $D$ and $D_{\,*}$ be domains in ${\Bbb B}^n/G$
and ${\Bbb B}^n/G_*,$ respectively, while $\overline{D}$ and
$\overline{D_{\,*}}$ are compacta.

\medskip
Let $\varphi\colon(0,\infty)\rightarrow (0,\infty)$ be a
nondecreasing function with a Calderon condition
\begin{equation}\label{eqOS3.0a}
\int_{1}^{\infty}\left[\frac{t}{\varphi(t)}\right]^
{\frac{1}{n-2}}dt<\infty\,,
\end{equation}
in addition, $K^{n-1}_O(p, f)$ is integrable in some neighborhood
$p_0\in \overline{D}$.
Then any homeomorphism $f\colon D\rightarrow {\Bbb B}^n/G_*$ in the
class $W^{1,\varphi}_{\rm loc}$ with a finite distortion
satisfies~(\ref{eq1F})--(\ref{eq*3!!}) at any point
$p_0\in\overline{D}$ for $Q(p):=c\cdot  K^{n-1}_O(p, f)$, where
$c>0$ is some absolute constant, and $K_O(p, f)$ is defined
in~(\ref{eq1B}). }
\end{proposition}

\medskip
Given Riemannian surfaces ${\Bbb S}$ and ${\Bbb S}_*,$ domains
$D\subset {\Bbb S}$ and $D_*\subset {\Bbb S}_*, $ a continuum
$E\subset D,$ $\delta>0$ and a measurable function $Q:{\Bbb
S}_*\rightarrow [0, \infty]$ we denote by $\frak{Sob}^Q_{E,
\delta}(D)$ the family of all homeomorphisms $f$ of $D$ onto $D_*$
such that $g=f^{\,-1}$ has a finite distortion, $K_O(p_*,
g)\leqslant Q(p_*)$ for almost all $p_*\in D_*,$ while
$\widetilde{h}_*(f(E))\geqslant \delta.$ The following statement
holds.

\medskip
\begin{theorem}\label{th3}
{ Let $D$ and $D_*$ be domains in ${\Bbb S}$ and ${\Bbb S}_*,$
respectively. Assume that, $D$ is a $2$-Loewner space,
$\overline{D_*}$ is a compactum in ${\Bbb S}_*,$ besides that, $D_*$
satisfies the condition~\textbf{A} and has at least two different
boundary points, and $Q\in FMO(\overline{D_*}).$ Let
$B_{\widetilde{h}}(p_0, \varepsilon_1)\subset D,$ where
$B_{\widetilde{h}}(p_0, \varepsilon_1)=\{p\in {\Bbb S}:
{\widetilde{h}}(p, p_0)<\varepsilon_1\}.$ Then there is $r_0>0,$
which does not depend on $f,$ such that
$$f(B_{\widetilde{h}}(p_0, \varepsilon_1))\supset
B_{\widetilde{h}_*}(f(p_0), r_0)$$
for all $f\in \frak{Sob}^Q_{E, \delta}(D, D_*),$ where
$B_{\widetilde{h}_*}(f(p_0), r_0)=\{p_*\in {\Bbb S}_*:
\widetilde{h}_*(p, f(p_0))<r_0\}.$ }
\end{theorem}

\medskip
\begin{proof}
The necessary conclusion follows from Theorem~\ref{th1} together
with Proposition~\ref{pr1}. Indeed, we may assume that
$\overline{B_{\widetilde{h}}(p_0, \varepsilon_1)}$ is a compactum in
$D.$ By Proposition~\ref{pr1A}, $D$ is locally connected and locally
compact metric space. By the same reason, ${\Bbb S}_*$ is locally
path connected and locally compact. It follows by
Proposition~\ref{pr3} that ${\Bbb S}_*$ is upper $2$-regular. In
addition, by the same Proposition
$$\widetilde{v}_*(B_{\widetilde{h}_*}(p^{\,*}_0, 2r))\leqslant C (2r)^2\leqslant
4C^2\widetilde{v}_*(B_{\widetilde{h}_*}(p^{\,*}_0, r))\,,$$
so that~(\ref{eq6}) trivially holds. Finally, any $f\in
\frak{Sob}^Q_{E, \delta}(D, D_*)$
satisfies~(\ref{eq2*A})--(\ref{eqA2}) with $Q=Q(p)$ by
Proposition~\ref{pr1}. The desired conclusion follows by
Theorem~\ref{th1}.
\end{proof}

\medskip
The corresponding analogue of Theorem~\ref{th3} may also be
formulated for the spatial case, see below. Given factor-spaces
${\Bbb B}^n/G$ and ${\Bbb B}^n/G_*,$ domains $D\subset {\Bbb B}^n/G$
and $D_*\subset {\Bbb B}^n/G_*, $ a continuum $E\subset D,$
$\delta>0,$ a nondecreasing function
$\varphi\colon(0,\infty)\rightarrow (0,\infty)$ and a measurable
function $Q:{\Bbb B}^n/G_*\rightarrow [0, \infty]$ we denote by
$\frak{OSob}^{Q, \varphi}_{E, \delta}(D)$ the family of all
homeomorphisms of $D$ onto $D_*$ such that $g=f^{\,-1}\in
W^{1,\varphi}_{\rm loc}(D_*)$ and $K^{n-1}_O(g, p_*)\leqslant
Q(p_*)$ for almost all $p_*\in D_*,$ while
$\widetilde{h}_*(f(E))\geqslant \delta.$ The following statement
holds.

\medskip
\begin{theorem}\label{th4}
{\, Let $n\geqslant 3$, and let $G_*$ be a some group of M\"{o}bius
automorphisms of the unit ball which acts discontinuously and have
no fixed points in ${\Bbb B}^n.$ Assume that, the group $G$, which
corresponds to the space ${\Bbb B}^n/G$, consists only from one
mapping $g(x)=x$.  Let $D$ and $D_{\,*}$ be domains in ${\Bbb
B}^n/G$ and ${\Bbb B}^n/G_*,$ respectively, while $\overline{D}$ and
$\overline{D_{\,*}}$ are compacta.

\medskip
Let $\varphi\colon(0,\infty)\rightarrow (0,\infty)$ be a
nondecreasing function with a Calderon condition~(\ref{eqOS3.0a}).
Assume that, $D$ is a $n$-Loewner space, besides that, $D_*$
satisfies the condition~\textbf{A} and has at least two different
boundary points, and $Q\in FMO(\overline{D_*}).$ Let
$B_{\widetilde{h}}(p_0, \varepsilon_1)\subset D,$ where
$B_{\widetilde{h}}(p_0, \varepsilon_1)=\{p\in {\Bbb B}^n/G:
{\widetilde{h}}(p, p_0)<\varepsilon_1\}.$ Then there is $r_0>0,$
which does not depend on $f,$ such that
$$f(B_{\widetilde{h}}(p_0, \varepsilon_1))\supset
B_{\widetilde{h}_*}(f(p_0), r_0)\qquad \forall\,\,f\in
\frak{OSob}^Q_{E, \delta}(D, D_*)\,,$$
where $B_{\widetilde{h}_*}(f(p_0), r_0)=\{p_*\in {\Bbb B}^n/G_*:
\widetilde{h}_*(p, f(p_0))<r_0\}.$ }
\end{theorem}

\medskip
\begin{proof}
The desired conclusion follows from Theorem~\ref{th1} together with
Proposition~\ref{pr2}. Indeed, we may assume that
$\overline{B_{\widetilde{h}}(p_0, \varepsilon_1)}$ is a compactum in
$D.$ By Proposition~\ref{pr1A}, $D$ is locally connected and locally
compact metric space. By the same reason, ${\Bbb B}^n/G_*$ is
locally path connected and locally compact. It follows by
Proposition~\ref{pr3} that ${\Bbb B}^n/G_*$ is upper $n$-regular. In
addition, by the same Proposition
$$\widetilde{v}_*(B_{\widetilde{h}_*}(p^{\,*}_0, 2r))\leqslant C (2r)^n\leqslant
2^nC^2\widetilde{v}_*(B_{\widetilde{h}_*}(p^{\,*}_0, r))\,,$$
so that~(\ref{eq6}) trivially holds. Finally, any $f\in
\frak{OSob}^{Q, \varphi}_{E, \delta}(D, D_*)$
satisfies~(\ref{eq2*A})--(\ref{eqA2}) with $Q=C\cdot Q(p_*)$ with
some absolute constant by Proposition~\ref{pr2}. The desired
conclusion follows by Theorem~\ref{th1}.
\end{proof}

\section{The version for Riemannian manifolds}

All necessary definitions related to Riemannian manifolds and
Orlicz-Sobolev classes may be found, for example, in~\cite{Sev$_2$}.
Given Riemannian manifolds ${\Bbb M}^n$ and ${\Bbb M}^n_*$ with
geodesic distances $d$ and $d^{\,*}$ and volume measures $v$ and
$v_*,$ respectively, domains $D\subset {\Bbb M}^n$ and $D_*\subset
{\Bbb M}_*^n, $ a continuum $E\subset D,$ $\delta>0,$ a
nondecreasing function $\varphi\colon(0,\infty)\rightarrow
(0,\infty)$ and a measurable function $Q:{\Bbb M}^n_*\rightarrow [0,
\infty]$ we denote by $\frak{OSob}^{Q, \varphi}_{E, \delta}(D)$ the
family of all homeomorphisms of $D$ onto $D_*$ such that
$g=f^{\,-1}\in W^{1,\varphi}_{\rm loc}(D_*)$ and $K^{n-1}_O(g,
p_*)\leqslant Q(p_*)$ for almost all $p_*\in D_*,$ while
$\widetilde{d}_*(f(E))\geqslant \delta.$ (Here $K_O(p_*, g)$ is
defined similarly to~(\ref{eq0.1.1}) in the corresponding charts,
i.e., $K_O(p_*, g)=K_O(\psi(p_*), F),$ where $F=\varphi\circ
g\circ\psi^{-\,1},$ $(U, \varphi)$ are local coordinates of $g(p_*)$
and $(V, \psi)$ are local coordinates of~$p_*$). The following
statement holds.

\medskip
\begin{theorem}\label{th3A}
{ Let $n\geqslant 3,$ and let $D$ and $D_*$ be domains in ${\Bbb
M}^n$ and ${\Bbb M}^n_*,$ while $\overline{D}$ and
$\overline{D_{\,*}}$ are compacta. Let
$\varphi\colon(0,\infty)\rightarrow (0,\infty)$ be a nondecreasing
function with a Calderon condition~(\ref{eqOS3.0a}). Assume that,
$D$ is a $n$-Loewner space, besides that, $D_*$ satisfies the
condition~\textbf{A} and has at least two different boundary points,
and $Q\in FMO(\overline{D_*}).$ Let $B_{d}(p_0,
\varepsilon_1)\subset D,$ where $B_{d}(p_0, \varepsilon_1)=\{p\in
{\Bbb M}^n: d(p, p_0)<\varepsilon_1\}.$ Then there is $r_0>0,$ which
does not depend on $f,$ such that
$$f(B_{d}(p_0, \varepsilon_1))\supset
B_{d_*}(f(p_0), r_0)\qquad \forall\,\,f\in \frak{OSob}^Q_{E,
\delta}(D, D_*)\,,$$
where $B_{d_*}(f(p_0), r_0)=\{p_*\in {\Bbb M}_*^n: d_*(p,
f(p_0))<r_0\}.$ }
\end{theorem}

\medskip
\begin{proof} The desired conclusion follows from Theorem~\ref{th1}.
Indeed, we may assume that $\overline{B_{d}(p_0, \varepsilon_1)}$ is
a compactum in $D.$ By the definition, $D$ is locally connected and
locally compact metric space. By the same reason, ${\Bbb M}_*^n$ is
locally path connected and locally compact. It follows by
\cite[Lemma~5.1]{ARS} ${\Bbb M}^n_*$ is upper $n$-regular, so
that~(\ref{eq6}) trivially holds. Finally, by~\cite[Theorems~3.1 and
4.2]{ARS}, any $f\in \frak{OSob}^{Q, \varphi}_{E, \delta}(D, D_*)$
satisfies~(\ref{eq2*A})--(\ref{eqA2}) with $Q=C\cdot Q(p_*)$ with
some absolute constant. The desired conclusion follows by
Theorem~\ref{th1}. \end{proof}

\medskip
{\bf Funding.} The work was supported by the National Research
Foundation of Ukraine (Project ``Analogues of Carath\'{e}odory and
Koebe-Bloch theorems for Orlycz-Sobolev classes'', Project number
2025.02/0010).

CONTACT INFORMATION

\medskip
{\bf \noindent Evgeny Sevost'yanov} \\
{\bf 1.} Zhytomyr Ivan Franko State University,  \\
Velyka Berdychivs'ka Str 40, 10 008  Zhytomyr, UKRAINE \\
{\bf 2.} Institute of Applied Mathematics and Mechanics\\
of NAS of Ukraine, \\
19 Henerala Batyuka Str., 84 116 Slov'yans'k,  UKRAINE\\
esevostyanov2009@gmail.com

\medskip
{\bf \noindent Valery Targonskii} \\
Zhytomyr Ivan Franko State University,  \\
Velyka Berdychivs'ka Str 40, 10 008  Zhytomyr, UKRAINE \\
w.targonsk@gmail.com

\medskip
{\bf \noindent Denys Romash} \\
Zhytomyr Ivan Franko State University,  \\
Velyka Berdychivs'ka Str 40, 10 008  Zhytomyr, UKRAINE \\
dromash@num8erz.eu

\medskip
{\bf \noindent Nataliya Ilkevych} \\
Zhytomyr Ivan Franko State University,  \\
Velyka Berdychivs'ka Str 40, 10 008  Zhytomyr, UKRAINE \\
Email: ilkevych1980@gmail.com


\begin{thebibliography}{99}


\bibitem{AFW} { T. \ Adamowicz, K. \ F\"{a}ssler,  {\rm and} B. \ Warhurst}:
\textit{A Koebe distortion theorem for quasiconformal mappings in
the Heisenberg group}, Ann. Mat. Pura Appl., {\bf 199} (2020),
147-186.

\bibitem{Af}
{E. S. \ Afanas'eva}: \textit{On boundary behavior of one class of
mappings in metric spaces,} Ukrainian Math. J., {\bf 66} (1) (2014),
16-29.

\bibitem{ARS}
{E. S. \ Afanasieva, V. I. \ Ryazanov {\rm and} R. R. Salimov}:
\textit{On mappings in the Orlicz-Sobolev classes on Riemannian
manifolds,} J. Math. Sci., {\bf 181} (1) (2012), 1-17.

\bibitem{AS}
{T. \ Adamowicz {\rm and} N. \ Shanmugalingam}:
\textit{Non-conformal Loewner type estimates for modulus of curve
families,} Ann. Acad. Sci. Fenn. Math., {\bf 35} (2010), 609-626.

\bibitem{CG}
{L. \ Carleson {\rm and} T. W. Gamelin}: \textit{Complex Dynamics,}
Springer, New York (1993).

\bibitem{Cr}
{M. \ Cristea}: \textit{On the limit mapping of a sequence of open,
discrete mappings satisfying an inverse Poletsky modular
inequality,} Complex Anal. Oper. Theory, {\bf 19} (2025), Art.~165.

\bibitem{FS}
{A. \ Franovskii {\rm and} E. \ Sevost'yanov}: \textit{On mappings
with the inverse Poletsky inequality in metric spaces,} Ukr. Math.
J., {\bf 166} (2) (2022), 453-479.

\bibitem{Fu}
{B. \ Fuglede}: \textit{Extremal length and functional completion,}
Acta Math., {\bf 98} (1957), 171-219.

\bibitem{He}
{J. \ Heinonen}: \textit{Lectures on Analysis on Metric Spaces,}
Springer, New York (2001).

\bibitem{Ku}
{K. \ Kuratowski}: \textit{Topology, Vol. 2,} Academic Press, New
York-London (1968).

\bibitem{MRV$_1$}
{O. \ Martio, S. \ Rickman {\rm and} J. \ V\"{a}is\"{a}l\"{a}}:
Definitions for quasiregular mappings, Ann. Acad. Sci. Fenn. Ser. A
I Math., {\bf 448} (1969), 1-40.

\bibitem{MRV$_2$}
{O. \ Martio, S. \ Rickman {\rm and} J. \ V\"{a}is\"{a}l\"{a}}:
\textit{Topological and metric properties of quasiregular mappings,}
Ann. Acad. Sci. Fenn. Ser. A I Math., {\bf 488} (1971), 1-31.

\bibitem{MRSY}
{O. \ Martio, V. \ Ryazanov, U. \ Srebro {\rm and} E. \ Yakubov}:
\textit{Moduli in Modern Mapping Theory}, Springer, New York (2009).

\bibitem{MS}
{O. \ Martio {\rm and} U. \ Srebro}: \textit{Automorphic
quasimeromorphic mappings in $\mathbb{R}^n$,} Acta Math., {\bf 135}
(1975), 221-247.

\bibitem{Ri}
{S. \ Rickman}: \textit{Quasiregular Mappings,} Springer, Berlin
(1993).

\bibitem{RV}
{V. \ Ryazanov {\rm and} S. Volkov}: \textit{On the boundary
behavior of mappings in the class $W^{1,1}_{\mathrm{loc}}$ on
Riemann surfaces,} Complex Anal. Oper. Theory, {\bf 11} (2017),
1503-1520.

\bibitem{Sev$_1$}
{E. \ Sevost'yanov}: \textit{On the local and boundary behavior of
mappings on factor-spaces,} arXiv:1905.06414.

\bibitem{Sev$_2$}
{E. \ Sevost'yanov}: \textit{On the local and boundary behaviour of
mappings of factor spaces,} Complex Var. Elliptic Equ., {\bf 67} (2)
(2022), 284-314.

\bibitem{SevSkv}
{E. \ Sevost'yanov {\rm and} S. \ Skvortsov}: \textit{On the
convergence of mappings in metric spaces with direct and inverse
modulus conditions,} Ukr. Math. J., {\bf 70}(7) (2018), 1097-1114.

\bibitem{SM}
{E. \ Sevost'yanov {\rm and} A. \ Markysh}: \textit{On
Sokhotski-Casorati-Weierstrass theorem on metric spaces,} Complex
Var. Elliptic Equ., {\bf 64}(12) (2019), 1973-1993.

\bibitem{ST$_1$}
{E. \ Sevost'yanov {\rm and} V. \ Targonskii}: \textit{On
convergence of homeomorphisms with inverse modulus inequality,} Eur.
J. Math., \textbf{11} (1) (2025), Art.~17.

\bibitem{ST$_2$}
{E. \ Sevost'yanov {\rm and} V. \ Targonskii}: \textit{An analogue
of Koebe's theorem and the openness of a limit map in one class,}
Anal. Math. Phys., {\bf 15}(3) (2025), Art.~59.

\bibitem{Skv}
{S. \ Skvortsov}: \textit{Local behavior of mappings of metric
spaces with branching,} J. Math. Sci., {\bf 254}(3) (2021), 425-438.

\bibitem{Vu}
{M. \ Vuorinen}: \textit{Conformal Geometry and Quasiregular
Mappings,} Lecture Notes in Math., vol. 1319, Springer, Berlin
(1988).

\end{thebibliography}
\end{document}